# Le coup d'œil du scarabée :
# Charles Nordmann, ou pourquoi la guerre est l'affaire des savants

David Aubin

Institut de mathématiques de Jussieu,
Université Pierre et Marie Curie
david.aubin@upmc.fr

Juillet 2010

VERSION NON DEFINITIVE, ne pas citer sans autorisation de l'auteur

## Introduction

15 septembre 1914 : le lectorat bourgeois de la *Revue des deux mondes* (RDM) a bien d'autres préoccupations que les débats intellectuels de salons feutrés qui le passionnent habituellement. La guerre a éclaté depuis un mois et demi. La bataille de Marne s'achève, la ligne de front se stabilise de l'Oise à l'Argonne, pour l'instant Paris est hors de danger… C'est alors que la RDM entrouvre ses pages aux horreurs des combats qui font rage. Le chroniqueur scientifique de la revue, l'astronome Charles Nordmann, y livre ses premières impressions de « combattant » [Nordmann 1914j]. A 33 ans, il est mobilisé dans le 7e bataillon du génie mais n'a pas encore vu le front. D'un point de vue qu'il veut « subjectif », mêlant avec aisance gouaille populaire et solennité patriotique, il raconte les heures graves que vit la France dans les semaines qui suivent la mobilisation. Dans les taxis et les trains qu'il emprunte, les villes et les villages qu'il traverse, Nordmann se fait le porte-parole de la résolution de tout un peuple, une résolution calme pimentée d'esprit gaulois. Il s'agit bien là du récit avec lequel on est familier : l'union sacrée qui gomme les doutes et les divisions, un engagement complet et spontané de tout un peuple qui ne vise qu'un seul but, la victoire sur l'ennemi héréditaire [Becker 1977].

Le cœur serré, Nordmann repense à l'Observatoire de Paris où il laisse ses travaux scientifiques en plan. L'astronome ne quitte pas sans mélancolie les étoiles et leurs secrets. Cherchant réconfort dans le vieil adage rabelaisien selon lequel science sans conscience ne serait que « trouble » de l'âme, il comprend que son statut de savant ne dépend que des circonstances sociales qui lui permettent d'exercer son métier. « L'univers avec ses merveilles, ses astres d'or, son immensité dans le temps et l'espace, ses transformations étonnantes, tout cela n'existe réellement que parce que cela est *pensé* par nous » [Nordmann 1914j, 131]. Plus que pour la France — ou parce la France en est pour lui la garante — c'est, « pour la pensée et ses droits que nous allons nous battre, pour la liberté de mieux étudier, mieux connaître, mieux utiliser la nature » [Nordmann 1914j, 132].

Vers le 15 août, dans le village d'Anjoutey, dernier poste avant l'Alsace allemande, c'est le choc. Pour la première fois, Nordmann croise des soldats qui reviennent du front, « quelques éclopés » :

> Ce sont presque des enfans. Ils ont l'air fatigué, mais une flamme dans les yeux. Ils se sont battus l'avant-veille […]. Notre canon de 75 a fait merveille. Des Prussiens qui s'avançaient en masse, en colonnes de quatre, ont été fauchés par lui soudain et sont

tombés… toujours en colonne de quatre. L'un de ces petits soldats nous raconte que, parmi les monceaux de cadavres prussiens, il en a vu un resté debout sur ses jambes, figé et comme pétrifié par la mélinite dans cette attitude [Nordmann 1914j, 144].

C'est, dans son récit comme dans la RDM, la première mention de la brutalité du combat. D'emblée, elle semble indissociable de ces fleurons français de la technique moderne que sont le canon de 75 mm et la poudre de mélinite. On perçoit l'impression désagréable que cela laisse au scientifique, à peine atténuée par le fait qu'il décrive l'anéantissement de soldats ennemis. La mort fauche la jeunesse des nations occidentales et son auxiliaire est la technique. Pour Nordmann, comme on le verra, il ne fait aucun doute que sous-tendant tout progrès technique, il y a ce qu'il appelle « la science », cette science pour laquelle il déclare aller au combat. La science, qui autrefois faisait la fierté de l'Europe moderne, est aujourd'hui complice du massacre : pire, c'est elle qui le rend possible !

Aujourd'hui, il y a consensus parmi les historiens au sujet de ce qu'on pourrait appeler la « culture de guerre » scientifique en France [Prochasson et Rassmussen 1996 ; Hanna 1996 ; Aubin et Bret 2003]. Selon ce point de vue, les savants français dans leur ensemble auraient, dès août 1914, spontanément offert le soutien de leur savoir et savoir-faire au gouvernement et aux autorités militaires. Suite aux rumeurs provenant des zones occupées et surtout en réaction au célèbre « appel au monde civilisé » des 93 intellectuels allemands, certains savants parmi les plus éminents se seraient alors rués dans une croisade intellectuelle en défense de la culture scientifique française. Puis, à partir de 1915, aurait enfin été organisée, sous l'égide des Painlevé et autres Borel, une véritable mobilisation scientifique coordonnée au niveau ministériel. Laboratoires, bureaux de calculs, usines, services techniques de l'armée auraient enfin été en mesure de coordonner leur action jusqu'à jouer un rôle important dans la mobilisation industrielle nécessaire à la victoire finale sur les Empires centraux.

L'investissement massif des savants dans l'effort de guerre, cependant, ne va pas sans profondément remettre en question les convictions d'autrefois. L'internationalisme de la science, le désintéressement des savants, la foi dans le progrès qu'elle incarne : tout cela est ébranlé par la férocité des combats. Pourquoi les savants offrent-ils si peu de résistance à leur enrôlement dans une guerre qui semble être la négation de leurs valeurs humanistes ? Pourquoi, au contraire, montrent-ils tant d'enthousiasme et de spontanéité dans l'offre de service ? Comment, plus généralement, s'accommode-t-on d'un monde qui contemple « avec stupeur cette régression qui mettait au service de la barbarie déchaînée tous les progrès scientifiques de la civilisation » [Mangin 1920, 775] ? Alors même que les communautés savantes consentent à travailler sur des techniques mortifères sous le contrôle direct de l'armée, on reconstruit une image publique de la science qui demeure globalement positive. Comment, donc, est-il possible que la Grande Guerre ait pu apparaître comme étant la première guerre scientifique sans que l'image de la science en tant que force de progrès n'en soit irrémédiablement ternie ?

Ainsi, alors que les communautés savantes consentent à travailler sur des techniques mortifères sous le contrôle de l'armée, il semble qu'on s'accommode relativement aisément d'une image publique de la science qui, en retreignant le jugement moral porté à son encontre aux applications qu'on en fait, reste globalement positive. Déjà présent avant-guerre, le discours selon lequel la science désintéressée — les mathématiques y compris — propose des innovations qui peuvent être utilisées à bon ou mauvais escient s'impose en quelques années et devient banal [Pestre 1984]. On nous l'a rappelé récemment : « Le consentement n'est jamais inné, il est toujours fabriqué » [Lambelet 2005, 298]. Depuis quelques années, de vives querelles ont éclaté parmi les historiens de la Grande Guerre. Il est significatif, à mon sens, que ces querelles se soient nouées autour de la question du consentement et de la contrainte. Elles opposent ceux qui font l'hypothèse que l'engagement dans le conflit est un acte volontaire, produit d'une « culture de guerre » [Audoin-Rouzeau et Becker 2000], à des représentants d'une « école de la contrainte », selon laquelle les participants de cette guerre en furent plus souvent victimes forcées[1] [Cazals 2002]. La manière dont scientifiques, militaires et profanes apprennent à vivre dans un monde où la science, autrefois force de progrès, peut être dévoyée est, elle aussi, fabriquée pendant la première guerre mondiale. Les questions se posent : les savants dans leur ensemble consentent-ils à cette mobilisation scientifique ou, au contraire, se fait-elle sous la contrainte des événements ? Et quel rôle joue l'image des sciences dans cette acceptation ?

Pour tenter d'apporter des réponses à ces questions, il est peu de trajectoires individuelles qui soient plus éclairantes que celle de Nordmann. Chez ce dernier, la mobilisation scientifique est en effet inséparable de sa propre expérience du combat, même brève, et des luttes idéologiques qui s'engagent à propos des sciences et de leur rôle dans la guerre moderne. En imaginant, dès l'automne 1914, un système de repérage des batteries ennemies par le bruit qu'elles émettent au moment de leur détonation, Nordmann est d'ailleurs l'un des tout premiers scientifiques qui s'engagent, autour d'un système technique, dans une collaboration active avec les autorités militaires. Ensuite, c'est un écrivain qui touche un public large et influent. Chroniqueur scientifique à la RDM depuis 1912 et pendant toute la durée de la guerre, Nordmann est bien placé pour observer, accompagner et, parfois, donner quelque impulsion à ce mouvement qui place les sciences au premier plan des débats intellectuels, politiques, industriels et militaires de l'époque. La seule interruption importante dans ses chroniques se situe entre septembre 1914 et août 1915, lorsqu'il est mobilisé dans l'artillerie, puis chargé d'effectuer des expériences de repérage par le son. Mais, comme on l'a vu, il n'abandonne pas son rôle de médiateur. Au contraire, il rédige ses « impressions d'un combattant » qui seront publiées en feuilleton dans la RDM, avant d'être réunies dans un livre publié pendant la guerre [Nordmann 1917a]. Plus qu'un point de vue individuel sur les opérations militaires, ces « notes de routes » seront pour lui l'occasion de mettre en évidence le caractère proprement scientifique de la conduite de la guerre en ce début du 20ᵉ siècle, non seulement dans les laboratoires de l'arrière mais aussi sur le front. En bref, Nordmann est un savant qui manie le canon, les rupteurs et la plume.

« Des circonstances curieuses et particulières front que j'ai été précisément à même d'éprouver tour à tour dans ce domaine [de l'invention] les sentiments du jugé et celui du juge… Les malfaiteurs deviennent, dit-on, les meilleurs policiers » [Nordmann 1916b, 692]. Acteur aussi bien que commentateur, Nordmann exprime donc, à chaud, une opinion sur la guerre et sur la science qu'il est intéressant d'analyser avec soin. Sauf en ce qui concerne le repérage par le son, l'astronome ne revendique jamais les premiers rôles. Tant dans ses chroniques que dans son témoignage, il jette, comme il dira dans la nécrologie de son maître Henri Poincaré, « le coup d'œil du scarabée sur le chêne superbe » [Nordmann 1912a, 332]. Modeste, ce regard est néanmoins précieux, car s'il est permis d'interroger la candeur de l'auteur, il ouvre une fenêtre sur la manière dont se construisent les discours sur la place de la science dans la guerre et la société moderne, discours qui trouveront une large résonance dans la société française d'après-guerre. Il nous renseigne sur la façon dont les anciennes conceptions épistémologiques, sociales et idéologiques sont ébranlées par l'expérience de la première guerre mondiale et les réajustements auxquels elles seront soumises. Par-dessus tout, il fait l'éclatante démonstration que la guerre est devenue l'affaire des savants.

---

[1] Notons que dans ces analyses la contrainte n'est pas toujours autoritaire et violente : il y a aussi des formes d'auto-contraintes comme les sentiments de camaraderie ou encore du fait de l'adhésion à certaines valeurs (patriotiques, par exemple).

**La démonstration de Nordmann : la guerre scientifique et ce qu'elle implique**

Parmi les témoignages de scientifiques et d'inventeurs sur la Première Guerre mondiale, il en est un qui, bien que terriblement acrimonieux, est particulièrement éclairant au sujet des difficultés liées à la pleine reconnaissance du rôle de ces derniers dans la guerre. Il s'agit d'un pamphlet au vitriol publié à compte d'auteur par Gorges Claude en 1919 intitulé *Politiciens et Polytechniciens* [Claude 1919]. Né en 1870, Claude est physicien, vulgarisateur et industriel : il a entre autres fondé la Société L'Air Liquide et développé l'industrie de l'éclairage au néon [Aubin 2003]. A la signature de l'Armistice, l'importance de sa contribution personnelle à l'effort de guerre, de même que celle des industries qu'il a fondées, est reconnue par son élection dans la toute nouvelle division de l'Académie des sciences dédiée aux applications des sciences à l'industrie. Pourtant, Claude digère mal d'avoir été incapable de convaincre les autorités des sections techniques militaires du bien-fondé de quelques unes de ses inventions, en particulier celles liés au repérage par le son. « Nos savants, écrit-il ! Ah! si vous saviez ce que l'on a fait d'eux ! […]Si vous saviez le peu qu'ils ont pesé pour une caste de gens monstrueusement préoccupés de leurs seuls privilèges, obstinément butés à écarter la concurrence possible » [Clause 1919, 33].

Son ressentiment n'est pas récent. Dès le 12 novembre 1915, le *Matin* de Paris publiait une lettre incendiaire sous le titre : « La guerre et nos savants : aux inventions barbares de la science allemande, il faut une riposte française ». Dans cette lettre, Claude ne cachait rien de l'amertume qu'il éprouvait vis-à-vis de « la modestie du rôle joué en ces jours tragiques par la science française ». Non que les savants français aient manqué d'enthousiasme patriotique ou de compétences techniques : on les a empêchés, affirmait-il, de remplir un devoir à la hauteur de leur talent. Claude pointe du doigt les services techniques de l'armée :

> Débordés — fatalement — par la multiplicité des problèmes scientifiques que les nécessités du front soulèvent chaque jour, on les aimerait, ces services, soucieux de reconnaître leur embarras et, pour mener à bien leur mission, cherchant parmi les compétences notoires des collaborateurs indispensables. […] Serait-ce trop de dire que chez ceux-là, nos savants, bien souvent, ont été accueillis sans bienveillance [Claude 1915] ?

Ainsi, selon Claude, aussi bien dans les récits rétrospectifs dans ses réactions « à chaud », la cause est entendue. A cause de leur incompétence technique et d'un esprit de corps malvenu, les autorités militaires n'ont pas accueilli comme elles auraient dû le faire les initiatives patriotiques des savants. Avaient-ils tort ? Parfois, il semble que les services militaires ou gouvernementaux aient su répondre aux exigences nouvelles plus efficacement que les scientifiques de l'université [voir Aubin 2000 et Aubin à paraître]. Mais, à la fin de la guerre, le constat qui s'impose est tout autre : c'est la recherche en mathématique, en physique, en chimie, bref en sciences fondamentales, qui entre désormais de plein pied dans les questions militaires et les affaires d'État. Président de l'Académie des sciences après avoir dirigé la France pendant la guerre, le mathématicien Paul Painlevé le souligne lors de la séance annuelle de remise des prix en décembre 1918 :

> Les mathématiques les plus abstraites ou les plus subtiles ont participé à la solution des problèmes de repérage et au calcul des tables de tir toutes nouvelles qui ont accru de 25 pour 100 l'efficacité de l'artillerie. […] [L]a guerre a rendu plus étroite et féconde l'union de l'Industrie et de la Science ; elle a mobilisé la Science au service de la Patrie [Painlevé 1919, 3-4].

L'expérience de Nordmann a-t-elle été bien différente de celle de Claude ? Dans l'état actuel de nos connaissances, on ne saurait le dire avec certitude. Il est par contre certain qu'à deux reprises Nordmann change de rôle. D'abord soldat, il consacre ensuite quelques mois à des recherches sur le repérage par le son. Puis, le 16 juillet 1915, il est officiellement affecté à l'arrière, en tant qu'adjoint au secrétaire général de la Commission supérieure des inventions intéressant la défense nationale[2]. Selon toute évidence, l'astronome ne cherchera plus à promouvoir ses propres inventions auprès des militaires, il jugera celles des autres. De cette plateforme privilégiée, il sera à même d'observer les mutations des rapports entre science et guerre.

Propagandiste avant tout, Nordmann cherche à faire une démonstration assez subtile, d'une rigueur toute mathématique si l'on veut, à propos de la guerre moderne et du rôle qu'y joue la science. Avec une belle régularité, il poursuit cet objectif sans jamais l'expliciter complètement. La lecture exhaustive de ses écrits sur ce sujet ne laisse pourtant aucun doute sur la constance remarquable avec laquelle il vise à prouver non seulement que la science est nécessaire à la guerre moderne, non seulement que les savants en sont maintenant les acteurs incontournables, mais d'abord et avant tout que la guerre est, dans tous ses aspects, devenue proprement scientifique. De là découle le fait que l'application des méthodes scientifiques à des fins militaires, loin de constituer un dévoiement des valeurs qu'on chérissait avant-guerre, en est l'aboutissement. Et c'est précisément par leur adhésion à ces valeurs positives, voire positivistes, beaucoup plus qu'à cause de leurs simples compétences techniques et leur inventivité d'ailleurs partagées, que les scientifiques acquièrent l'autorité qui leur confère toute légitimité de se prononcer sur tous les aspects de cette guerre.

Afin de mieux saisir cette démonstration, nous procéderons par étapes. Dans un premier temps, nous utiliserons les chroniques scientifiques que Nordmann a écrites avant guerre afin de faire apparaître les valeurs qu'il associe alors à « la science ». Puis, nous verrons comment il s'y prend, dans le témoignage qu'il rapporte du front, pour démontrer le caractère scientifique de la Grande Guerre. Nous reviendrons ensuite sur son travail de recherche à propos du repérage, montrant, avec lui, que si les principes scientifiques jouent un rôle de premier plan dans ce conflit, ce n'est peut-être pas par ses compétences techniques que le savant est le plus utile dans cette guerre. Les chroniques scientifiques publiées par Nordmann pendant la première guerre mondiale et dans l'immédiat après-guerre nous permettront d'affiner cette conclusion. S'il en arrive à penser que dans la société qui s'annonce le savant prend place entre le général et l'homme d'État, s'il exprime clairement que la science ne peut servir de base à la morale, il n'abandonne pas pour autant le rêve d'une science poétique, désintéressée et internationale.

**Les chroniques d'avant la guerre : la science, l'internationalisme, les applications et le sentiment religieux**

---

[2] Historique de la Commission [des inventions] : « Par note de service en date du 16 juillet 1915 n° 1920BL, le Ministre de la Guerre désigne le sous-lieutenant Nordmann, astronome à l'observatoire de Paris, est affecté au 12e régiment d'artillerie, pour être adjoint au Lieutenant-colonel, secrétaire général de la Commission des Inventions. » (AN F17 17486). L'auteur remercie Christian Gilain pour cette citation. Voir aussi le dossier militaire de Nordmann conservé au SHD-Terre 6-Ye-52777. Sur l'histoire de la Commission des inventions, voir [Roussel 1989] et [Anizan 2006].

Astronome adjoint à l'Observatoire de Paris, Nordmann a 33 ans lorsque le conflit éclate[3]. En 1903, il a soutenu, devant Poincaré, une thèse de doctorat ès sciences physiques. Deux ans plus tôt, âgé d'à peine 20 ans, il a gravi le Mont Blanc dans le but de détecter des ondes radio émises par le soleil. Ses observations à la station des Grands-Mulets à 3 100 mètres d'altitudes avaient été un échec, qu'il suppose dû à l'absorption des ondes par les hautes couches de l'atmosphère. Ses travaux scientifiques ne lui apportent aucune gloire et il semble s'être détourné assez vite de la recherche active et n'avoir pas été des plus assidus à l'Observatoire [Saint-Martin 2008 ; Aubin & Le Gars, à paraître]. Suite à la célèbre nécrologie qu'il consacre à son « maître » Poincaré [Nordmann 1912a], il devient chroniqueur scientifique à la RDM. Auparavant, il y avait déjà publié un article sur la vie et la mort des étoiles [Nordmann 1910], et était l'auteur de quelques articles de vulgarisation ou de nécrologies dans la *Revue scientifique* et dans le quotidien *Le Matin*. Le 1[er] décembre 1912, dans un contexte de crise économique à la RDM, et « pour répondre au passionnant intérêt du public au mouvement scientifique » [procès-verbal de l'assemblée des actionnaires, cité in Karakatsoulis 1995, 45] il reprend le flambeau du chroniqueur scientifique laissé par le biologiste Albert Dastre : « un honneur dont je sens le péril » [Nordmann 1912b, 697].

La rhétorique du vulgarisateur est convenue : « le premier devoir du chroniqueur scientifique est d'oublier complètement sa personnalité » [Nordmann 1914j, 129]. C'est pourquoi il est nécessaire d'être précis dans l'analyse du discours qu'on en fait. Comme son « parrain » Dastre, l'astronome embrasse le rôle du traducteur : si la science procure à l'initié « la plus haute des jouissances intellectuelles », les profanes y ont difficilement accès. Pourtant il veut, comme son prédécesseur, être non seulement l'« interprète » mais aussi la « personnification même de la science » [Nordmann 1917m, 661]. Il ne s'adresse pas qu'au profane, mais aussi aux professionnels, aux « rats de laboratoire », dont il veut exprimer les « idées plus ou moins inconscientes […] sous une forme explicite » [citation non référencée de Dastre, in Nordmann, 1917m, 667]. C'est donc en tant que fabricant d'une opinion de la communauté savante, une opinion plus ou moins consensuelle et qui resterait accessible aux non-initiés, que Nordmann se pose d'emblée.

En 1912, la RDM a perdu son lustre d'antan, mais elle n'en reste pas moins un carrefour intellectuel des plus actifs. Depuis sa fondation en 1829, la RDM porte une attention soutenue aux enjeux scientifiques, aux côtés d'autres préoccupations où les questions militaires restent bien présentes. Toujours libérale, elle est cependant devenue plutôt modérée, catholique et patriote [Broglie 1979, 15]. Au moment où Nordmann la rejoint, le caractère de cette « vieille institution littéraire » peut être résumé par la caractérisation qu'en fait un critique littéraire en 1916 : « une sorte d'académie. […] Elle est essentiellement traditionnelle, et elle admet et se permet toujours quelques traits de modernité ; elle est partisan d'un ordre bien établi, et d'une autorité ferme, et elle est toujours, avec cela, très libérale ; l'esprit chrétien y subsiste et la liberté et la diversité des croyances y sont admises et y sont respectée[4] ». Depuis quelques dizaines d'années déjà, la RDM est au premier plan des débats entre science et religion [Loué 1999 ; Loué 2003]. Si elle ne s'oppose jamais frontalement à l'image progressiste de la science, la rédaction de la RDM reste circonspecte vis-à-vis des dérives du positivisme et réaffirme souvent qu'il serait illusoire de vouloir fonder une nouvelle morale sur les sciences modernes. Les chroniques de Nordmann prennent ainsi place dans un forum intellectuel proche de milieux catholiques qui redoutent particulièrement l'empiètement du scientifique sur la morale [Fontana 1990 ; McMillan 1995]. Comme on le verra, les ressorts autrefois mobilisés pour fustiger les dérives scientistes qui s'attaqueraient aux domaines réservés de la religion, de la culture classique et de la morale seront remis au goût du jour et serviront maintenant à dénoncer les pratiques allemandes.

---

[3] Fils d'un fabricant d'horloge originaire de Besançon, Charles Nordmann est né le 18 mai 1882 à Saint-Imier dans le canton de Berne. Il mourra en 1940. Cf. [Cru 1929, 453] et son dossier militaire au SHD.
[4] Ces extraits proviennent d'un article d'Émile Faguet publié dans le *Gaulois* au début de 1916 [cité par Broglie 1979, 351-352].

Les vingt-quatre articles que Nordmann écrit pour la RDM jusqu'en septembre 1914 couvrent un large spectre de domaines scientifiques. Bien qu'astronome, il ne se limite pas au domaine déjà large couvert par le personnel des observatoires. Il rapporte les principales nouvelles dans le domaine de l'astronomie, la chimie et la physique, la biologie et la médecine, l'océanographie et l'exploration de la planète, etc. Les applications techniques ou médicales des sciences y figurent en bonne place, tout comme les conventions internationales et les célébrations importantes. Chacune de ses chroniques faisant au moins douze pages, l'astronome prend le temps d'introduire ses sujets de façon ample en les agrémentant souvent d'un long résumé historique.

Nordmann n'est pas mathématicien et les chroniques qu'il publiera dans la RDM ne traiteront que rarement et obliquement des mathématiques. Cela s'explique d'autant mieux que les mathématiques sont, pour lui comme pour François Arago près d'un siècle auparavant, l'un des principaux obstacles à la large compréhension du progrès scientifique. C'est peut-être même l'« appareil mathématique terriblement rébarbatif » qui en bloque l'accès au plus grand nombre. Il faut donc traduire : « il faut faire apercevoir par-dessus le squelette dont [les scientifiques] se contentent, la figure réelle et harmonieuse de la vérité » [Nordmann 1910, 535]. Il lui arrive même d'écrire qu'il ne s'agit guère que d'un « appareil pédantesque » [Nordmann 1923, 8]. Certes, il admet que les mathématiques sont parfois indispensables : il est « difficile de faire sur l'aviation, sans le secours des mathématiques, un exposé d'une diaphanéité vraiment aérienne » [Nordmann 1916f, 699]. Mais souvent, il préfère s'en passer imaginant accroître le nombre de ses lecteurs. S'il est d'accord que pour parler correctement de balistique, il faut « le secours austère de quelques équations », il veut « renoncer pourtant ici à cet appareil ésotérique, et […] n'exposer ces choses […] que dans le simple langage de tout le monde » [Nordmann 1918g, 460]. Dans les années 1920, cette attitude fera de lui l'un des grands vulgarisateurs de la relativité : « si Einstein est un trésor, un horrible troupeau de reptiles mathématiques en éloigne le curieux. Qu'il y ait en eux, comme dans certaines gargouilles gothiques, une secrète beauté, c'est certain. Mais il vaut mieux, armés du fouet éclatant qu'est le verbe, les chasser loin de nous et accéder jusqu'aux splendeurs einsteiniennes par le clair et nobles escalier du langage français[5] ».

Mais ce n'est pas uniquement parce qu'elles se prêtent mal à la vulgarisation que Nordmann exclut les mathématiques. Répétant souvent la phrase de Poincaré — « l'expérience est la source unique de la vérité » [Poincaré 1900, ***] — il n'est qu'admiration pour les mathématiciens la mettent en application [Nordmann 1916b, 689]. Dans sa hiérarchie de valeurs, les mathématiques ne sauraient servir de guide : « dans le domaine de la science, le fait est souverain, et la pensée, si vive qu'elle soit, ne peut généralement que se traîner dans son sillage impérieux » [Nordmann 1915e, 215]. Pourtant, Nordmann ne peut manquer de s'apercevoir que les mathématiques sont florissantes en France et que c'est dans ce domaine que l'apport des jeunes générations paraît le plus brillant. Selon lui, ce choix est fait par défaut. Les moyens en France sont souvent disponibles pour les scientifiques en fin de carrière, non pour les jeunes : « Seules les sciences mathématiques échappent à cette règle, et c'est pourquoi sans doute ces sciences sont les seules où nous soyons, non pas seulement les égaux des premiers, mais hors pair » [Nordmann 1917j, 936n].

---

[5] [Nordmann 1921a], p. 10. Grand défenseur de la clarté du verbe, il écrit ailleurs : « les Fontenelle, les d'Alembert, les Arago nous ont enseigné qu'il n'est arcane dans les sciences où le net langage français ne sache porter à l'occasion un clair rayon de lumière » [Nordmann 1918g, 460].

Même s'il écarte les mathématiques de ses chroniques, Nordmann n'en demeure pas moins convaincu que les sciences les plus avancées sont celles qui sont les plus mathématisées. C'est dans la nécrologie de Poincaré qu'apparaît le plus clairement le point de vue de Nordmann sur les sciences. Il s'attache d'abord à y décrire l'homme : c'est l'image idéale du savant qu'il construit. Poincaré donne « une impression singulière de spiritualité et d'impérieuse douceur » [Nordmann 1912a, 332]. On affirme qu'il est peu au fait des applications auxquelles ses travaux mathématiques peuvent donner lieu : le nécrologue prétend que les anecdotes véhiculées à ce sujet trahissent plutôt un trait essentiel de son appréhension de « la Science pour l'Art » [Nordmann 1912a, 336]. Le Poincaré de Nordmann est un génie « sensible à la beauté sous toutes ses formes, pourvu qu'elles fussent nobles » : musique, peinture, poésie, mais aussi les sciences qui lui procurent des « ravissemens esthétiques » [Nordmann 1912a, 336]. De fait, il est lui-même un artiste, « le Michel-Ange de la pensée » qui ne voit pas les détails, qui ne s'attarde pas aux minces récoltes, mais « dans les recoins les plus obscurs, les plus inabordables des choses, qui savait d'emblée, à grand coups de ciseau, faire jaillir des allées nouvelles pleines de lumières et de fleurs » [Nordmann 1912a, 337]. Il résume l'homme ainsi :

> Son âme est avant tout celle d'un artiste et d'un poète. Ces vues si profondes et si vraies vont un peu à l'encontre des idées classiques sur le type du "savant", respectable, certes, mais un peu caricatural, avec son cerveau mécanique, son œil que les lunettes traditionnelles rendent aveugle à toute beauté et son cœur où la nature a déposé, au lieu de sensibilité, une table de logarithme à sept décimales [Nordmann 1912a, 343].

C'est à ce sentiment exacerbé pour l'esthétique, à sa grande faculté de généralisation qu'il attribue le fondement du génie mathématique de Poincaré. Si son œuvre en astronomie est présentée comme « gigantesque », c'est sur sa philosophie, qui a « profondément remué tout ce qui pense dans le monde » [Nordmann 1912a, 350], que s'attarde l'auteur. C'est chose remarquable que Poincaré parvienne par sa pensée à « toucher une époque réaliste et vulgaire où les conflits des intérêts priment plus que jamais sur ceux des idées » [Nordmann 1912a, 350]. Cette époque, Nordmann la décrit souvent comme étant traversée par des conflits profonds au sujet de la valeur de la science. Il met en scène deux camps qui s'opposent sur la question des rapports entre la science et le monde réel, entre la science et la morale ; ces oppositions lui seront bien utiles lorsqu'il s'agira de rendre compte de la manière dont les sciences interviennent dans la conduite de la guerre. Entre deux extrêmes, le chroniqueur occupe le juste milieu[6].

A propos de la science et du monde, il met en scène, d'un côté, « le matérialisme, le rationalisme, le positivisme et le scientisme » [Nordmann 1912a, 351] qui prétendent parfois « nous fournir une image définitive en quelque sorte dogmatique de l'Univers » [Nordmann 1912a, 356], et, de l'autre, une vision « pragmatique » selon laquelle la science n'aurait de valeur que comme moyen d'action et ne serait qu'une « création artificielle, incertaine, contingente et qui ne nous apprend rien sur la réalité objective » [Nordmann 1912a, 359]. Partant de la philosophie poincaréenne, il adopte une position très proche de celle qui sera plus tard celle des divers courants structuralistes, à savoir « que la seule réalité objective, ce sont les rapports, les relations des choses », à propos desquels la science « nous enseigne tout ce qu'il y a en eux de purement objectif » [Nordmann 1912a, 360 ; Aubin 1997]. On ne peut donc pas s'appuyer sur Poincaré pour fonder une conception « agnostique » de la science, bien que la foi qu'on peut avoir en telle ou telle théorie soit limitée. Poincaré offre « des raisons de douter, mais en même temps d'agir, et d'aimer le beau et le vrai, bien qu'ils soient peut-être inaccessible » [Nordmann 1912a, 352].

---

[6] « *A priori* d'ailleurs, — et il en est ainsi hélas ! dans toutes les choses humaines, — il est certain que ni l'une ni l'autre [de deux positions contradictoires] ne doit être parfaitement bonne ni parfaitement mauvaise, que la plus juste doit contenir quelque part contestable et que dans la plus fausse doit surnager un brin de vérité » [Nordmann 1918h, 934].

Cette position épistémologique rejoint aussi sa conception des rapports entre la science et ses applications. Peu développée dans l'article sur Poincaré, cette question traverse son mode d'écriture. Proche des vulgarisateurs du 19e siècle comme Louis Figuier, Nordmann mélange allègrement science la plus pure et applications les plus concrètes. Apparaît souvent le sentiment d'une grande continuité d'objet, d'intention et de méthode entre les divers intervenants, fussent-ils savants, médecins ou industrialistes. Souvent, ses articles présentent d'abord l'histoire, puis les théories les plus récentes sur un thème scientifique avant de s'attarder à leurs conséquences médicales ou industrielles. Les deux paires de chroniques qu'il consacre au froid [Nordmann 1913a ; 1913b] et aux ultraviolets [Nordmann 1913k ; 1913l] sont construites sur ce modèle. Dans « Les tendance et les progrès récens de la chimie », il écrit :

> On a beaucoup discuté pour savoir si les sciences sont nées de nécessités pratiques ou de ce besoin désintéressé de savoir qui tourmente les hommes. Pour les sciences mathématiques, la question pouvait se poser, et il est encore des esprits ingénieux que n'a point convaincus la démonstration que fit Henri Poincaré de leur origine utilitaire. Pour la chimie, il n'y a point de doute possible : non seulement elle est née des besoins matériels de l'humanité, mais elle puise encore en eux sa principale raison d'être [Nordmann 1913c, 218].

L'utilitarisme n'est pas en soi condamnable. Mais là n'est pas la principale raison des sciences. Dans les années d'avant-guerre, le thème qui prédomine dans les chroniques nordmanniennes a en effet trait à la morale et au sentiment religieux qu'inspire l'activité scientifique. Conservateur et catholique, le public de la RDM s'attend à ce que son chroniqueur scientifique, même d'origine juive [on trouve mention de Nordmann dans Gygès 1965, 220], soit investi d'une mission. Son rôle est de mettre en évidence le caractère moral et spirituel de la science, ou, plus précisément, à mieux en circonscrire le périmètre. Encore une fois, soulignons que les thèmes qui apparaissent ici fourniront les ressources indispensables permettant de penser le clash entre science et guerre dans les années qui vont suivre.

C'est à nouveau dans la nécrologie de Poincaré que ces questions sont les mieux. Nordmann pose explicitement la question : « la vérité scientifique peut-elle être en conflit avec la morale ? » Dans ce débat, le chroniqueur adopte la position d'arbitre entre ceux qui voient dans la science « une école d'immoralité » et d'autres qui, à l'instar d'Anatole France, prétendent voir poindre une « Morale issue des sciences naturelles » [Nordmann 1912a, 361-362]. Sa position est sans équivoque : « la science ne peut servir de base à une morale impérative » [Nordmann 1912a, 363]. Cette constatation n'implique portant pas que la science soit condamnée à la neutralité morale : au contraire, la sélection naturelle peut, par exemple, nous aider à comprendre le sentiment du juste et de l'injuste, l'équilibre entre l'égoïsme et l'altruisme qui caractérise l'« honnête homme » [Nordmann 1912a, 364]. De fait, il affirme, comme Poincaré, qu'il n'y a pas de vérité qui soit dangereuse pour la société et que « le meilleur remède contre une demi-science, c'est plus de sciences » [Poincaré, cité par Nordmann 1912a, 364].

L'astronome de l'Observatoire s'oppose farouchement à l'idée que le savant doive abandonner tout sentiment religieux. Au contraire, l'homme de science qui a gardé la foi doit être envié, car il est plus complet que l'idéaliste ou que le réaliste. A l'instar du géologue Albert de Lapparent[7], le chroniqueur milite donc, non pour une démarcation nette entre les domaines relevant du scientifique et du religieux, mais pour une conciliation. Il n'est pas question de « "compartimenter" [l'univers] en zones d'influences dont les unes pourraient être divines et les autres non » [Nordmann 1913f, 708]. Au contraire : « la science nous montre, en effet que l'Univers est un tout ordonné, cohérent, harmonieux ; et c'est par là plus encore que par ses dimensions qu'il est grandiose ; c'est par là qu'il est mystérieux et divin » [Nordmann 1912a, 367].

Ce sentiment poétique et religieux du cosmos affleure souvent dans ses chroniques. La célébration des 25 ans de l'Institut Pasteur offre l'occasion de revenir sur ce thème [Nordmann 1914b]. Rappelant que son fondateur était chrétien et croyant, citant ses écrits sur le non empiètement de la science sur la religion, il réaffirme que la métaphysique ne peut se réduire à la physique. Dans le conflit qui oppose alors matérialistes et néo-vitalistes à propos de savoir si le vivant est soumis à d'autres forces que la matière inanimée, le chroniqueur assume une nouvelle fois la position médiane [Aubin 2008]. Dans le doute, il est impossible d'être affirmatif, dit-il, mais du même souffle, fait-il remarquer, on a pas besoin de la position vitaliste pour affirmer sa foi dans le divin, que la loi de la gravitation universelle n'est pas moins mystérieuse que la vie et la mort des êtres humains et que « ce qu'il y a de divin dans l'Univers, c'est précisément que son harmonie soit gouvernée par des lois immuables, universelles et ne souffrant pas d'exception » [Nordmann 1913f, 708].

Enfin un dernier thème apparaît dans les chroniques d'avant-guerre, qui aura aussi son importance dans l'interprétation de l'expérience de la guerre : l'internationalisme et la place de la France dans le paysage mondial des sciences. Dans ces chroniques, les savants français occupent, on le comprend aisément, une place centrale, mais il est à noter qu'avant 1914, leurs travaux sont souvent replacés dans un contexte international. Des domaines de recherche où la France est peu présente, comme la conquête des pôles, sont même largement commentés [Nordmann 1913d ; 1913e]. Il souligne aussi le rôle joué par Darboux dans l'instauration d'instances scientifiques internationales [Nordmann 1914a].

Si la science se déploie dans un contexte international, Nordmann n'est pas insensible au rang que son pays y occupe. A l'occasion de la 5e conférence internationale des poids et mesures en 1913, il fait le point sur la situation internationale du système métrique, « clair reflet du génie français » [Nordmann 1913l, 697]. Voilà qui « inflige un démenti sans réplique aux esprits chagrins qui sans cesse et sans raison geignent sur l'amoindrissement de notre influence dans le monde » [Nordmann 1913l, 683]. Conquête de l'esprit, le système métrique s'est montré plus fort que les baïonnettes et règne maintenant sur de nouvelles aires géographiques. C'est à cette occasion, rappelant l'histoire des savants de la Convention, que le chroniqueur cite une phrase attribuée à Pasteur, dont il aura amplement l'occasion de se resservir dans les années qui suivent : « si la science n'a pas de patrie, le savant doit en avoir une » [Nordmann 1913l, 683].

---

[7] Voir *Science et apologétique : conférences faites à l'Institut catholique de Paris, mai-juin 1905*, 2e éd., Paris : Bloud, 1905 ; cité in Nordmann, « Poincaré », 367.

Lorsque le prix Nobel est attribué aux chimistes français Paul Sabatier (1854–1941) et Victor Grignard (1871–1935) en 1912, Nordmann en profite pour contrer l'impression répandue selon laquelle « la chimie est aujourd'hui une science exclusivement allemande » [Nordmann 1913c, 217]. La science est bien « une et indivisible », mais, précise-t-il, « ce noble internationalisme ne doit pas nous empêcher de d'observer que les apports scientifiques des peuples reflètent à leur manière les qualités et les défauts de chacun d'entre eux ». Bien sûr, les découvertes des sciences exactes restes « impersonnelles », il n'empêche : on peut y retrouver « l'âme » de ceux qui les ont faites. Dans ce but, Nordmann entreprend donc, plus d'un an avant le début de la première guerre mondiale, de tracer le portrait de ce que la chimie moderne doit au caractère français de certains de ses pratiquants. Les savants français y apportent sans cesse des idées nouvelles, tout en « laissant à d'autres le soin et l'honneur de les appliquer et d'en tirer les lointaines conséquences ». D'où le fait que l'industrie chimique française ne soit pas même la deuxième au niveau mondial, « bien que la France ait produit ces dernières décades la plus riche moisson de découvertes chimiques qu'on puisse imaginer » [Nordmann 1913c, 218]. On sait à quel point ces arguments deviendront monnaie courante dans la propagande produite ensuite par les Pierre Duhem, les Émile Picard et autres Edmond Perrier. En 1913, Nordmann ne semble guère troublé par cet état de choses : au contraire, après avoir insisté sur la synthèse de produits nouveaux, il insiste plutôt sur la « poésie grandiose et mystique [qu'il y a] dans le laboratoire du chimiste » [Nordmann 1913c, 228]. Dans l'imaginaire qu'il dépeint, le savant français, même chimiste, reste plus près de la poésie que des usines.

Tous ces thèmes bien développés avant-guerre serviront, au moment du déclenchement des hostilités, de ressources indispensable à la construction de nouveaux discours. Elles trouveront une nouvelle vie pendant et après la Grande Guerre. La croisade intellectuelle française assimilera les deux extrêmes que sont l'abstraction mathématiques et l'utilitarisme pur à la science allemande, forcément dévoyée. On pourra d'un même souffle célébrer les succès français dans l'application de principes scientifiques à la guerre et dénoncer les « barbaries » allemandes, puisque les premiers restent, contrairement aux secondes, encadrés par une morale d'inspiration religieuse, sinon chrétienne. Peu à peu, l'expérience du combat et celle de la recherche militaire feront naître de nouvelles façons de concilier désintéressement et utilitarisme, spiritualité et pessimisme, internationalisme et guerre mondiale. Voyons d'abord comment Nordmann transformera son expérience du combat en argument polémique.

### Témoin : la guerre comme « expérience scientifique »

« Le phénomène bataille est comme tous les phénomènes naturels justiciable de l'expérimentation et de la critique scientifique ». Cette phrase de Nordmann est la toute première des nombreuses épigraphes que Norton Cru place en tête du célèbre livre qu'il consacre aux « témoins » de la Grande Guerre [Cru 1929]. Dans l'esprit de Cru, il s'agit d'appliquer la méthode scientifique à l'étude de l'histoire de la guerre. Très sévère dans son appréciation des témoignages qui exagèrent ou, pire, répandent des rumeurs, Cru place le récit de Nordmann parmi ceux dont il recommande vivement la lecture. Il est, souligne-t-il, l'un des deux seuls docteurs ès sciences à avoir livré le récit de leur expérience dans les tranchées. Sous la plume de Cru, il s'agit-là d'une qualité rare qui permet aux auteurs de reconnaître les limites du témoignage personnel, à savoir son caractère subjectif, mais aussi sa grande force qui est de livrer un récit vrai parce que basé sur une expérience vécue. L'autre docteur, Georges Kimpflin, y insiste : « Dans la guerre présente, [la] vérité est aussi celle des combattants. [...] Le combattant a des vues courtes [...]. Mais parce que ses vues sont étroites, elles sont précises ; parce qu'elles sont bornées, elles sont nettes. Il ne voit pas grand'chose, mais il voit bien ce qu'il voit. Parce que ses yeux et non ceux des autres le renseignent, il voit ce qui est » [Kimpflin 1920, 12-14 ; cité par Cru 1929, p. 353].

« Mis en demeure de rédiger des impressions » par le directeur de la *Revue des deux mondes* [selon Cru 1929, 182], Nordmann produira aussi un document s'appuyant essentiellement sur ses propres observations. A l'instar du pigeon de Jean de La Fontaine, il affirme : « J'étais là, telle chose m'advint ». Reporter de guerre, il veut changer de ton. Il abandonne le style impersonnel du chroniqueur scientifique et promet de livrer ses « impressions, c'est-à-dire les reflets dans une âme particulière d'événemens particuliers ». Ses impressions, écrit-il encore, seront, « à l'encontre des écrits scientifiques, lentement "subjective"[8] ». Mais il ne se départit pas de sa mission pédagogique : « Je parlerai, néanmoins, le moins possible de moi, juste ce qu'il faut pour qu'on n'oublie point que ces heures ont été vécues, que ces choses ont été vues, senties » [Nordmann 1914j, 129-130].

Pour le lecteur d'aujourd'hui, pourtant, le tableau que Nordmann nous dresse de la guerre peut néanmoins sembler convenu et peu personnel, sauf en ce qui concerne des détails peu significatifs. Le pédagogue prend souvent le dessus sur le témoin. Au fur et à mesure de ses chroniques, la tâche qu'il s'est fixée se précise. Suite au premier volet publié, on l'a vu, dès le 15 septembre 1914, les articles suivant paraissent le 1er novembre 1914 [Nordmann 1914j ; 1914k], les 15 mars, 15 juin, et 15 juillet 1915 [Nordmann 1915a ; 1915b ; 1915c], et enfin le 15 novembre 1916 [Nordmann 1916l]. Légèrement remaniées et réarrangées, ces chroniques paraissent dans un ouvrage préfacé par le général Robert Nivelle [Nordmann 1917a]. Il semble que l'ouvrage ait connu un certain succès puisque selon l'exemplaire microfilmé à la BNF, toujours daté de 1917, il en serait déjà à sa 9e impression.

Dans leur ensemble, ces différents écrits, nous l'avons dit, font l'effet d'une démonstration. La guerre est scientifique : dans *Le Matin,* il va jusqu'à écrire qu'elle est « une vaste expérience [...] d'astronomie » [Nordmann 1917n] ! Ou plus facétieusement : « La guerre est une équation à plusieurs inconnues qui admet des "racines" [...] multiples. Ce qui importe surtout, — tous ceux qui ont mordu un peu au biberon des mathématiques me comprendront, — c'est d'écarter les "solutions imaginaires" » [Nordmann 1918h, 938]. Mais contrairement aux autres chantres de la guerre scientifique en 1914 (H.G. Wells, Salomon Reinach, Jules Violle, par exemple), Nordmann ne se contente pas d'énoncer cet avènement : il le montre comme une réalité du champ de bataille. Il ne s'agit pas seulement d'affirmer que les sciences jouent désormais un rôle de premier plan, non seulement dans le développement de nouvelles armes et la mise au point de nouvelles thérapeutiques, mais au cœur même de l'action militaire. C'est ce que souligne le général Nivelle dans sa lettre-préface, datée du 1er décembre 1916. Louant le talent et le « pittoresque » du style de Nordmann, Nivelle espère qu'ils aideront les Français à comprendre « cette longue et formidable guerre dont la complexité exige le concours sans limites et sans réserve de toutes les sciences, de tous les arts, de toutes les industries, de toutes les énergies matérielles et morales ». Le général redoute cependant certains excès de l'esprit scientifique de l'auteur : il serait, écrit-il, illusoire de « chercher dans une formule qui ne serait applicable qu'au procédé de guerre actuel, qui dispenserait de réfléchir et de vouloir, la solution des problèmes variés à l'infini que pose la guerre » [Nivelle, in Nordmann 1917a, i-v].

---

[8] Il ajoute : « Qu'on me pardonne ce mot, que j'emprunte aux Allemands, mais nous avons bien d'autres choses à leur emprunter, et même à leur prendre... » [Nordmann 1914j, 130]

Les premières « impressions » que livre Nordmann sont celles d'un flâneur de la guerre. Soldat, il est resté libre. Fort de ses capacités équestres et de sa connaissance de l'allemand, le voilà cavalier. Le front l'attire « ne serait-ce que par curiosité et amour de l'imprévu » [Nordmann 1914j, 138]. L'amusement sera de courte durée : « nous nous sommes adaptés avec une rapidité prodigieuse aux péripéties de la grande tragédie qui a succédé à nos petites comédies d'hier » [Nordmann 1914k, 18]. Dès le deuxième article daté du 1ᵉʳ novembre 1914, il adopte un ton plus « scientifique », le récit est moins linéaire et laisse plus de place à l'analyse. Il est vrai qu'à ce moment, il a déjà quitté le front et rédige sans doute d'après ses notes. Déjà, le soldat se laisse aller à des réflexions sur l'histoire, sur la météorologie, et même sur une éclipse de soleil qui, le 21 août, est cachée par les nuages. Son troisième article, publié le 1ᵉʳ mars 1915, laisse entrevoir que son auteur a entretemps eu d'autres occupations : « Les choses sont si étrangement organisées sur cette planète que la pensée, du moins la pensée spéculative, et l'action sont presque exclusives l'une de l'autre » [Nordmann 1915a, 319]. Les événements qu'il a vécus depuis sont, affirme-t-il, « tout à fait romanesque » et il promet d'en faire un jour l'histoire qui « tiendra à la fois d'un conte de Voltaire et d'un roman de Jules Verne, — avec, hélas ! le talent en moins » [Nordmann 1915a, 319]. On regrette, évidemment, que Nordmann n'ait pas — ou si peu — tenu promesse[9] !

C'est à ce moment qu'est trouvé le ton qui sera celui d'*A coups de canons*. Ses prochains articles ne seront pratiquement pas retouchés lors de la publication du livre en 1917. S'il n'admet pas encore, comme il le fera à ce moment-là, que « la guerre n'est pas le poème romantique que certains avaient rêvé, mais une sombre prose réaliste, [et que] les valeurs ont été renversées » [Nordmann 1917a, 9], il enfourche dès lors un nouveau cheval de bataille :

> regardons froidement en face la dure réalité. Assez de sublimes enfantillages. […] La guerre a cessé d'être un art pour être une science expérimentale comme la physique. Regardons-la d'un œil achromatique. […] Tâchons de « comprendre » le mécanisme de cette guerre [Nordmann 1917a, 9-10].

Dès lors, il reprend son rôle de chroniqueur scientifique, mais sans l'exaltation qui le caractérisait auparavant. Il va rechercher « non le "pourquoi" mais le "comment" de cette guerre, non sa sombre poésie, mais son mécanisme, comme on ferait d'une réaction chimique » [Nordmann 1917a, vii]. Puisque « l'artillerie est bien une science » [Nordmann 1915a, 320 = 1917a, 11-12], il s'attache à décrire par le menu les opérations d'artillerie : le défilement, le pointage, l'observation des cibles, le téléphone, les instruments d'optique, les calculs balistiques, les obus que l'on envoie et ceux qu'on reçoit. Jusqu'aux observatoires militaires qui deviennent comparables à ceux qu'il fréquentait auparavant :

> C'est étonnant comme le nombre des observatoires s'est multiplié depuis quelques temps sur le territoire, et surtout tout le long de cette mince ligne qu'on appelle le front, sur laquelle déferlent nos énergies et qui profile l'armure de la France. Ce n'est point dû à un renouveau soudain des études astronomiques : ce ne sont point les étoiles qui en sont la cause, mais des êtres qui n'ont, hélas ! pas grand'chose de céleste ni d'éthéré : messires les Boches ; car il y a observatoire et observatoire, comme il y a fagot et fagot, comme il y a vérité et affirmation germanique.

> Donc il n'est point aujourd'hui de groupe d'artillerie, et généralement même de batterie, qui n'ait son observatoire d'où l'on règle, observe et dirige les coups. La consigne y est la même que dans les temples à coupoles où naguère nous disions la messe aux étoiles : "Faire des observations ou en recevoir." D'autre part, les télescopes sont-ils autre chose que des canons idéalisés ? Aussi la transition a-t-elle été toute naturelle et le voyage très rapide qui nous fit descendre soudain de Sirius sur ce coin de bonne terre franque, face aux Teutons [Nordmann 1915b, 888-889 = 1917a, 86-87].

---

[9] Il réitérera d'ailleurs cette promesse à plusieurs reprises [Nordmann 1916i, 228 ; 1917f, 945 ; 1919g, 702], mais ne s'exécutera jamais sauf par bribes [Nordmann 1916i ; 1916j ; 1920a ; 1928].

Si les opérations d'artillerie apparaissent dès lors comme des activités scientifiques à part entière, le travail scientifique lui-même, ses résultats et applications diverses deviennent aussi partie prenante de la guerre : « ce n'est pas le moindre paradoxe de cette guerre que de voir la chimie, cette chimie dont nos ennemis étaient si fiers […] nous donner sur eux un avantage décisif. Belle matière à philosopher sur la science en général, et la science allemande en particulier, sur leur rôle dans l'art de s'entre-massacrer et leur influence sur le bonheur de l'espèce humaine » [Nordmann 1915a, 336 = 1917a, 47]. Il ne voit pas sans tristesse le thème de la « faillite de la science » qui ressurgit dans les journaux. Signalons que c'est dans la *Revue des deux mondes*, à côté des chroniques de Nordmann, que paraissent, à cette époque, d'autres attaques bien connues contre la science allemande [Duhem 1915 ; Picard 1915]. La réponse à cette critique est vite trouvée : ce n'est pas la science qui a failli, pas même la science des Allemands, mais la manière dont ils l'ont faite ! « [R]ien ne permet de solidariser la science avec leurs sophismes sauvages » [**ref ???**].

Quant aux Français, ils ne doivent pas trop maudire « les moyens de destruction perfectionnés que la science a mis aux mains des chevaliers de la "Kultur" » puisque ce sont les mêmes moyens qu'ils ont à disposition pour se défendre : « ce sont eux précisément qui nous permettent de dominer dans l'artillerie » [Nordmann 1915a, 337]. En définitive, il n'a que faire des polémistes de l'arrière : « il est puéril et, en ce moment, malfaisant de ranimer de vieilles polémiques d'amphithéâtres pendant que les soldats se battent. La science n'est ni morale ni immorale ; on l'a démontré cent fois » [Nordmann 1915a, 338 = 1917a, 50]. Ailleurs, il affirmera ironiquement que le combattant se moque bien de ces querelles de privilégiés « dont l'héroïsme est surtout typographique » [Nordmann 1916a, 457]. Au contraire, c'est l'heure de l'union sacrée de la science et de la morale religieuse pour la défense de la France et de ses idéaux :

> Que celui qui porte dans son âme un idéal religieux, que celui qui porte l'amour ardent de la science, que celui plus heureux qu'illumine l'une et l'autre de ces torches intérieures, que chacun d'eux puissent croire, sans l'angoisse du doute, qu'en se battant pour la France, il se bat aussi pour l'étoile idéale qui guide ses pensées [Nordmann 1915a, 338 = 1917a, 51].

Le chroniqueur se bat pour la France, car c'est elle qui incarne le mieux la conception quasi-religieuse de la science qu'il défend. Il se donne pour mission d'utiliser sa plume afin de galvaniser la résolution de son lectorat. Son but ultime, c'est la Victoire qui « est fille du peuple, musclée et fière : la poudre aux yeux ne l'impressionne que si elle est pyroxylée ». Le moyen de l'atteindre s'exprime brièvement : « Travaillons ». C'est la conclusion de son livre [Nordmann 1916l, 414 = 1917a, 244]. Voyons maintenant comment le travail du savant peut, selon Nordmann, contribuer à la victoire.

**La technique des savants : le repérage par le son**

Le repérage des batteries ennemies par le son est l'un des grands succès militaires de la recherche scientifique française pendant la première guerre mondiale. Avant même que le conflit ne touche à sa fin, il fait, dans de nombreux cercles, figure d'exemple [par exemple Burgess 1917]. On loue ce cas de coopération réussie entre savants et militaires non seulement parce qu'il débouche sur des actions jugées efficaces mais surtout parce qu'il indique les façons concrètes dont les scientifiques peuvent mettre leurs compétences techniques au service de l'armée. Car, s'il est vrai que la première guerre mondiale représente une opportunité nouvelle pour la science [Millikan 1919], elle l'est d'autant plus pour les savants. Avant même que la victoire ne soit acquise, les scientifiques des nations alliées peuvent prétendre y avoir contribué de manière essentielle. Cette prétention, il faut le souligner, se fait au détriment des autres acteurs de l'innovation militaire, les inventeurs et les officiers des services techniques [Roussel 1989 ; Aubin 2003 ; Galvez-Béhar 2005 ; Aubin à paraître].

Nordmann est sans conteste l'un des tout premiers à avoir imaginé le principe du repérage par le son et à avoir conduit des expériences pour transformer ce principe d'acoustique en des opérations susceptible d'être effectuée à grande échelle dans des conditions de combat. Au début d'octobre 1914, un épais brouillard recouvre le front d'Alsace : incapable de repérer ses cibles, l'artillerie est contrainte de rester inactive. Nivelle soumet le problème à ses officiers. Simple canonnier, Nordmann lui aurait immédiatement fourni la solution. Enthousiaste, Nivelle écrit : « Si cette solution théorique devait recevoir une réalisation pratique, cela constituerait une telle supériorité qu'il convient de tout mettre en œuvre pour aboutir. Et il faut aboutir rapidement[10] ».

A la base du repérage par le son, il y a un principe qui est mathématiquement simple. On cherche à mesurer et à exploiter le léger décalage de temps dans la perception du son par plusieurs observateurs placés en des endroits différents dans le but de déterminer, par triangulation sonore, la position d'une batterie ennemie en train de tirer. Comme l'explique le colonel Emmanuel Vallier (1849–1921), « la pièce se trouve sur une hyperbole ayant pour foyer les positions des deux observateurs[11] » [Vallier 1915, 298]. Etant donné trois observateurs, on peut tracer deux hyperboles : « tout mathématicien déterminera sans peine la position de la pièce qui a tiré. Solution théorique facile, mais qu'il s'agissait de transporter dans la pratique » [Arthur-Lévy 1923, 438].

Envoyé à Paris le 20 octobre 1914 afin d'y conduire des expériences, Nordmann s'installe à l'Institut Marey au bois de Boulogne où il collabore avec son directeur Lucien Bull (1876–1972). Le 17 novembre, on procède à un essai grandeur nature dans la banlieue ouest de Paris. Testant deux systèmes différents l'un basé sur l'écoute humaine, l'autre utilisant un galvanomètre, devant une commission présidée par le général Joseph Gallieni, commandant du camp retranché de Paris à laquelle participe également Painlevé, Nordmann arrive à déterminer en une heure environ la position des pièces ayant tirés deux salves de six coups chacun, avec une précision de 20 à 40 mètres. Trois semaines plus tard, on est prêt à tester le dispositif sur le front. Une première batterie allemande dont la position était jusque là inconnue est détectée, le 8 décembre 1914, sur le front de l'Aisne, par la « section de repérage Nordmann (S.R.N. n° 1). » [Nordmann 1928, 29]. On utilise alors une méthode peu perfectionnée où les observateurs actionnent un bouton lorsqu'ils perçoivent le son du canon transmettant électriquement leur observation au poste de contrôle. Cette méthode, on s'en doute, reste peu précise, sauf si les observateurs sont assez éloignés les uns des autres [Claude 1917, 98-101]. Nivelle organise aussitôt les sections de repérage. Dès le 20 décembre 1914, des officiers sont formés aux nouvelles méthodes par Nordmann.

---

[10] Cette lettre est citée dans le rapport officiel sur les expériences concernant le repérage d'une batterie par la méthode de M. Nordmann (18 novembre 1914), signé par les généraux Clergerie, Desaleux, Miquel-Dalton, le colonel Girodon, et M. Paul Painlevé [reproduit *in extenso* in in Nordmann 1920a ; 1928, **pages\*\*\***]. Le journal de campagne du 5e RAC commandé par Nivelle (SHD-Terre, 26N913, 2) mentionne le brouillard mais ne révèle rien à propos du repérage.
[11] On peut à juste titre s'étonner que Vallier, par ailleurs secrétaire de la Commission des Inventions ait ainsi « defloré avec une superbe inconscience » ce procédé dans une publication aussi accessible : « Un an après, les Allemands multiplieront les sections de repérage. Coïncidence ? Je ne sais pas. Mais ce que je sais bien, c'est qui fut advenu si le colonel Vallier, né en un autre temps, avait rencontré sur sa route quelque membre du Comité du Salut Public ! » [Claude 1919, 95].

Entre-temps, « avec une noble émulation » [Nordmann 1928, 30], divers établissements scientifiques et techniques s'empressent de développer des systèmes concurrents, tous basés sur les mêmes principes. Quatre systèmes sont développés durant la guerre : celui du Service géographique de l'Armée (appelé système T.M. puisque développé par la télégraphie militaire), le système Dufour (professeur de physique à l'École normale supérieure), le système Cotton-Weiss, et le système Bull[12]. On sait aussi que d'autres scientifiques de grand renom y ont contribué, parmi eux : Ernest Esclangon [1925], Émile Borel, Henri Lebesgue, etc. Selon une note émanant du G.Q.G., datée du 21 décembre 1915, les différentes méthodes auraient donné des résultats à peu près comparables [citée in Nordmann 1928, 30]. Adapté par Lawrence Bragg, le procédé Nordmann est employé par les Anglais jusqu'à la fin de la guerre [Hartcup 1988, 68-76]. En 1918, l'armée française comptera 45 sections de repérage par le son dont l'effectif total s'élèvera à peu près à 6000 hommes : « des services scientifiques [qui] exigent pour fonctionner […] un grand nombre d'officiers compétents aidé d'un personnel subalterne spécialisé » [Bourgeois 1920, 684].

Comme on le verra dans le chapitre consacré au général Léon Bourgeois, qui dirige le Service géographique de l'Armée, la mise en opération des sections de repérage par le son constitue l'un des grands succès des sections techniques militaires [Schiavon *** ; Bourgeois 1920 ; Arthur-Lévy 1923]. Le problème du repérage est complexe, car il implique l'articulation d'un grand nombre d'éléments disparates. Une fois le principe acquis, les détails physiques du phénomène (en particulier, l'onde de choc) et les conséquences des perturbations météorologiques ou des configurations topographiques doivent être analysées avec soin, de délicats appareils électromécaniques doivent être développés, des procédures complexes doivent être établies afin de coordonner l'action des observateurs, des cartographes et des canonniers. Si pour certaines de ces tâches les scientifiques universitaires semblent incontournables, il en est d'autres que les inventeurs industriels et les ingénieurs militaires peuvent aussi bien, et parfois mieux, remplir.

C'est pourquoi l'histoire du repérage par le son constitue, dès le départ, un enjeu social que se disputent les divers acteurs de l'invention : mathématiciens et physiciens, ingénieurs militaires et inventeurs civils. Après-guerre, les récits qu'en font les acteurs laissent transparaître les frictions qui n'ont pas manqué d'éclater pendant la guerre. Au conservatoire des arts et métiers, le 22 février 1920, Bourgeois rend compte des succès du SGF : reconnaissant l'intuition de Nordmann, il s'évertue pourtant à ne pas le nommer [Bourgeois 1920]. C'est avant l'Armistice que l'astronome, quant à lui, a cherché à ce que sa double priorité soit reconnue : le premier, il aurait repéré des batteries tirant à blanc ; le premier des batteries ennemies[13]. Le 1er avril 1917, Nivelle lui envoyait une lettre confirmant sa priorité, que Painlevé avait déjà reconnue lors de la séance du 9 février 1916 de la Commission des inventions : « c'est M. Nordmann qui a, le premier, mis la chose sur pied ». Après enquête approfondie confiée à Charles Maurain (futur directeur de l'Institut de physique du globe), le directeur des Invention Jules Breton écrit à Clemenceau, ministre de la Guerre, en vue de la nomination de Nordmann au titre d'officier de la Légion d'honneur. Mais si Nordmann tient à ce que sa contribution technique soit pleinement reconnue, il attache plus d'importance à montrer que le rôle des savants dans cette guerre dépasse largement le cadre de la technique.

---

[12] Ces différents système, à l'exception du sytème Dufour, sont décrits dans une notice du Service géographique de l'Armée (sections de repérage par le son), *Repérage et réglage du tir par le son : description sommaire des divers systèmes employées dans l'Armée française en mai 1917*, Paris : Impr. nationale, 1917 [BNF cote 8° V-Pièces 28492]. Un document ronéotypé décrivant les caractéristiques techniques du système Cotton-Weiss, daté du 25 avril 1917, est aussi disponible à la BNF, sous le titre « Repérage des batteries par le son : système Cotton-Weiss », cote f° V 7153. Voir aussi [Schiavon 2003 ; sa thèse*** ; son chapitre].

**Un chroniqueur scientifique en guerre**

« Tout en devenant une question de matériel, la guerre demeure comme tout ce qui est de l'homme, un problème cérébral » [Nordmann 1917a 244]. Problème cérébral, la guerre sera donc le domaine de ceux qui cultivent la pensée. D'après la conception idéale que Nordmann se fait de la science, on comprend que les scientifiques soient par conséquent appelés à y prendre une place centrale. La RDM, « observatoire privilégié pour étudier la culture pendant la guerre » [Loué 2005, 184], est l'endroit tout désigné pour compléter sa démonstration et faire adopter ses conclusions par un publique large et influent. Membre du noyau dur rédactionnel, Nordmann est le troisième auteur en termes du nombre d'articles publiés entre 1914 et 1919. La guerre offre l'occasion à la RDM de renouveler en profondeur sa rhétorique et d'étendre considérablement son lectorat qui s'accroît, avec plus d'un million d'exemplaires vendus en 1919, de 75 % au cours du conflit [Karakatsoulis 1995]. Du fait de la place qu'y occupe le débat moral autour de la science, on mesure l'importance que revêtent les chroniques que rédige Nordmann.

Au début du conflit, la RDM est incertaine sur l'attitude qu'il convient adopter. Puis, le 15 octobre 1914, paraît un article d'Émile Boutroux, intitulé « Les Allemands et la guerre », dans lequel est introduit la problématique des rapports entre la science et la morale, qui sont profondément affectés par la guerre : « La brutalité est ici calculée et systématisée : […] il y a une barbarie savante » [Boutroux 1914, **page***]. On assiste alors au retour des valeurs anciennes de la RDM, le combat contre la prétention scientifique à intervenir dans le débat intellectuel, le combat « pour la responsabilité de l'intellectuel encadrée par une morale, à savoir la religion catholique, au-dessus de sa liberté, et *in fine*, contre l'idée que la science puisse constituer le fondement d'une morale » [Loué 2005, 184]. Si la guerre permet de faire coïncider nation et religion dans l'image d'une France éternelle, elle le fait au détriment de l'image de la science : c'est le vieux thème de la faillite de science qu'on ressuscite.

Cette conclusion, on l'a vu, déplaît profondément à Nordmann. Amorale plutôt qu'immorale et tout à fait compatible avec la morale religieuse sera donc la science de Normann, y compris cette science de guerre qu'il présente sans relâche à ses lecteurs, chaque mois, d'août 1915 à la fin du conflit : « ce que l'homme a su arracher pour la guerre à l'impassibilité tour à tour féconde ou meurtrière de la science » [Nordmann 1916k, 218].

---

[13] En 1923, on affirme bien que Nordmann fut le premier. Cependant, on indique que l'astronome Ernest Esclangon envoie un mémoire au SGA (qui se trouve alors à Bordeaux), dès le 20 septembre 1914. Dans la même semaine, M. Driencourt, ingénieur- hydrographe de la marine, et le colonel Ferrié, directeur de la télégraphie sans fil s'occupent de la question et la soumettent au SGA. Dans le même temps, Painlevé transmet à Bourgeois les résultats des essais fait à Paris par Nordmann. « L'idée parcourut-elle le monde scientifique ou naquit-elle spontanément dans les cerveaux qui tous n'aspiraient qu'à seconder les efforts de la défense nationale ? » [Arthur-Lévy 1293, 439].

> N'est-elle pas étrange et suggestive cette utilisation belliqueuse de ces ondes invisibles et muettes [les ondes radio de la TSF] qui ceignent la terre de leur houle silencieuse, et qu'on nous avait présentées, naguère, lorsque nous écoutions les douces ânonneries des utopistes, comme ne pouvant et ne devant jamais servir qu'à des fins pacifiques ? La vérité, c'est qu'il n'est aucune des conquêtes de la science […] qui ne puisse […] servir […] au besoin de combattre. Il n'y a pas de bons ou de mauvais outils ; il n'y a que de bons ou de mauvais travaux, et surtout de bons ou de mauvais artisans » [Nordmann 1919f, 938].

Puisque la science est indépendante de la morale, il lui faut un cadre. Pour Nordmann, ce cadre ne peut que celui qui est fourni par les valeurs de la nation. Les valeurs qu'il associait autrefois à la science en général deviennent des valeurs bien françaises : celles qui justement vont servir de garde-fous contre les dérives du scientisme utilitarisme et matérialiste.

Au début de 1916, un renouveau éditorial se fait sentir à la RDM. A la mort du directeur de la revue, c'est René Doumic, catholique qui s'affiche comme tel dès l'École normale, hostile au socialisme, qui en prend les commandes. Sans abandonner son conservatisme politique, la RDM s'ouvre à certaines idées nouvelles et assume une « nouvelle vocation » [Broglie 1979, 348]. On assiste alors à une certaine « littérisation de la guerre » et le « nationalisme littéraire » remplace le catholicisme comme idéologie de base [Loué 2005, 185 ; voir aussi Karakatsoulis 1995, 47ff]. Comme Doumic déclare à ses actionnaires :

> Appuyée sur sa tradition la *Revue* peut aller sans crainte vers les nouveautés et s'en assimiler ce qu'elles ont de meilleur. Faire une *Revue* qui, fidèle à son passé, devienne plus que jamais vivante, actuelle, variée, ce programme […] est celui qu'impose la situation [Karakatsoulis 1995, 63].

Pendant les heures sombres de la guerre mondiale, la RDM assume donc la mission de guide spirituel pour une nouvelle génération de décideurs politiques et militaires. Comme en témoigne le maréchal Foch dix ans plus tard : « Devant ces horizons chargés de brumes, nous avons trouvé le précieux concours de la *Revue des deux mondes*. Qu'il s'agît de hausser les cœurs ou d'éclairer les esprits, de remonter le moral ou d'asseoir et de fondre des opinions divergentes d'apparence, c'est un magnifique rayon de chaleur et de lumière qu'elle nous apportait » [cité in *Livre du centenaire* 411]. Dans la RDM où les questions militaires étaient auparavant déjà souvent soulevées et analysées avec finesse, furent publiés de nombreuses études générales sur l'armée, son organisation, le haut commandement, l'armement, et la condition des combattants. Les récits abondent couvrant tous les grands événements et tous les aspects du combat. De nombreux témoins y rapportent leur expérience du front.

Le 1er août 1915, Nordmann retrouve donc son poste de chroniqueur scientifique. Pendant toute la durée du conflit, il publiera près de quarante articles mensuels exclusivement consacrés à la science ou à la technique de guerre. Sur ces 40 études, 5 sont consacrées aux problèmes d'agriculture et d'alimentation, 8 à l'artillerie (technique et matériel), 4 à l'aviation ou l'aérostation, 5 à la chimie, 8 à la médecine, et pas plus de 5 ne s'occupent de questions politiques ou culturelles[14]. Si le thème principal de sa chronique a radicalement changé, entièrement tourné vers la lutte armée qui fait rage, les conceptions qui sous-tendent les chroniques de Nordmann ne témoignent guère que d'un infléchissement par rapport à celles d'avant-guerre. Le souci du vulgarisateur ne disparaît toutefois pas[15] : au contraire l'intérêt du public est plus vif que jamais. Depuis la comète de Halley, huit ans auparavant, il n'y a pas eu de « phénomène relevant de la science et comparable, par l'intérêt et les discussions qu'il a soulevés, au bombardement prodigieux de Paris par une pièce de longue portée » [Nordmann 1918d, 935].

---

[14] Voir la liste complète dans la bibliographie.

Le travail entrepris dans son témoignage de « combattant » se poursuit dans les chroniques de la RDM. Quoiqu'il déplore la censure et l'autocensure qui l'empêche trop souvent de dévoiler toute l'ampleur de la mobilisation scientifique, Nordmann donne à voir, à un public assez large, les multiples recoins de cette guerre où se cache la science[16]. Il n'est point, écrit-il par exemple, « d'engin qui enferme en ses flancs [autant] de science diaboliquement pervertie que la torpille » [Nordmann 1916k, 217]. La science est bien omniprésente dans le conflit actuel :

> Il n'est pas douteux que la monstrueuse Bellone [déesse de la guerre] qui préside aux tueries d'où sortira un monde rénové et affranchi, a pris aujourd'hui les traits de la science. Science, l'invention et l'extraction des explosifs qui propulsent les projectiles et en multiplient la puissance en mille éclats meurtriers ; science, la fabrication des aciers spéciaux d'où jaillissent les obusiers, les canons, les fusils, les mitrailleuses ; tous les engins mortels ; science, la transmission instantanée par le téléphone, le télégraphe avec ou sans fil, des ordres et des renseignements, qui lie et scelle au cerveau des chefs tous les membres épars de l'immense armée en campagne ; science, tous les procédés de topographie, de cartographie, de télémétrie, de repérage, qui font découvrir les objectifs cachés et permettent de les détruire ; science, l'emploi des avions de guerre, munis de tous les perfectionnements qui les lient à travers l'espace à leur point d'attache et assurent leur marche et la justesse de leur tir ; science, le sauvetage opératoire et le traitement de nos blessés et de nos malades qui chaque jour restitue de nouveaux bataillons à la France [Nordmann 1915h, 701].

La « domination monstrueuse des machines dans la lutte actuelle » [Nordmann 1916k, 217] ont accru l'importance de la technique. Certes la science elle-même en profite, certains phénomène d'illusion acoustique inconnus des physiciens ne l'étant pas des balisticiens : « Cela prouve que la science pure a toujours à gagner à ne point perdre contact avec les sciences appliquées » [Nordmann 1916j, 701]. De même à propos de l'aviation, il réfute l'idée « qu'elle avait pris la science au dépourvu, qu'elle n'était que le fruit de l'empirisme et que les savans n'étaient arrivés à la rescousse avec leurs *x* qu'après la bataille […]. Rien n'est plus faux » [Nordmann 1916f, 697]. Mais, c'est bien « la science » qui, dans la guerre scientifique de Nordmann, a le premier rôle, de même que c'est le savant qui comprend vraiment ce qui est en train de se passer.

La tentation serait forte, concède-t-il, d'y voir une perversion de l'entreprise qu'on glorifiait avant guerre. La civilisation ne serait qu'une « barbarie confortable » [Nordmann 1915h, 699-

---

[15] L'article sur l'agriculture rappelle ce qu'est l'azote et son importance pour la croissance des végétaux [Nordmann 1916c], celui sur l'avancement de l'heure montre la part jouée par les observatoires dans la mesure du temps [1916e], celui sur l'aviation explique le jeu des poussées qui fait voler l'avion [1916f], ses articles sur la poudre sont rigoureusement pédagogiques et historiques [1915d et 1915e]…
[16] A travers ces chroniques l'occasion s'offre d'ailleurs l'occasion d'explorer plus à fond la question de la censure et de l'autocensure en matière d'innovation technique [voir par exemple Nordmann 1915d, 684 ; 1915f, 691 ; 1916b, 691 ; 1916e, 217 ; 1918a, 445 ; 1918d, 940-941 ; 1919f, 927 ; 1918g, 458 1918h, 935].

700]. Vraie déesse de cette guerre, la chimie ne s'est-elle pas pervertie ? « C'est elle qui propulse dans les airs l'acier coupant et le cuivre lourd. […] La chimie étend aujourd'hui d'une façon sinistre […] sa mission qui est […] de décomposer les corps : les corps bruts des la nature minérale comme aussi les beaux corps si souples des jeunes guerriers » [Nordmann 1915f, 696]. Non, quand on parle de la guerre scientifique, il vaut mieux faire abstraction de ses convictions intimes. « Il vaut mieux en effet, si on veut l'observer avec intérêt, considérer ce phénomène du point de vue de la physique que de la morale » [Nordmann 1917g, 457]. Pourtant, tout espoir n'est pas perdu. En s'intéressant à la chirurgie de guerre, à la vaccination contre la typhoïde, à partir de 1917, c'est-à-dire aux manières nous pas de tuer, mais de combattre la mortalité, il retrouve parfois des accents lyriques : des découvertes « simples et lumineuses, qui sont nées de la guerre, comme, dans un, orage atroce et qui voile la douceur bleue du ciel, on voit jaillir parfois des éclairs » [Nordmann 1917m, 946].

Attaché à la Commission des inventions, Nordmann tire une première leçon de la guerre scientifique : « il faut que les prêtres de la science quitte eux aussi leur tour d'ivoire pour voler au tocsin » [Nordmann 1916b, 697]. Car ces ressources permettant d'organiser l'invention existent, ce sont les infrastructures scientifiques du pays. On doit utiliser les laboratoires des facultés pour seconder à ceux des services techniques de l'armée qui sont débordés, en particulier par la demande de nouveaux matériels. Mais ce ne sera pas uniquement pour leurs compétences techniques qu'on fera appel aux savants et à leurs infrastructures. Puisque les inventions non sollicitées qui arrivent spontanément sur les bureaux la Commission sont le plus souvent inutiles, il est nécessaire d'organiser l'invention. Il faut que les scientifiques œuvrent en collaboration avec, mais à côtés des, structures militaires, car ce n'est pas toujours le cas qu'on « apprécie davantage ce qui est inclus sous le képi que ce qui est cousu dessus » [Nordmann 1916b, 695]. Le principe de base, « c'est la meilleure utilisation des compétences », car « les supérieurs ne sont pas toujours les gens supérieurs » [Nordmann 1916b, 694]. Il faut ajoute-t-il non sans ironie « faire en sorte que le gallon, qui est une mesure de capacité chez les Anglais, en soit une aussi de l'autre côté de la Manche » [Nordmann 1916b, 695]. Les scientifiques rendront les militaires compétents !

Dès décembre 1915, il fournit de fait une analyse très clairvoyante de l'évolution des rapports entre science et État que cela induit. A l'instauration du ministère des inventions, il écrit ainsi : « C'est la première fois, à ma connaissance, que dans un document gouvernemental on admet officiellement la science à jouer un rôle dans les affaires de l'État […]. L'institution d'un organisme national destiné exclusivement à faire participer la science aux nécessités de l'heure, n'est donc rien moins qu'une sorte de révolution » [Nordmann 1915h, 698]. Il y a « paradoxe », bien sûr, à vouloir organiser les inventions [Nordmann 1916b, 687]. Le seul rôle à laquelle ce ministère peut prétendre, c'est « seulement à les utiliser, à les ordonner, à les adapter aux circonstances, à les dégager du réseau de barbelé des formalités administratives où se déchirent parfois leurs ailes délicates, à les défendre contre la routine, la bêtise, l'envie, le plagiat, en un mot à les traduire de concepts en faits » [Nordmann 1916b, 688]. Bref, il faut libérer les scientifiques des contraintes qui les empêchent d'agir aussi efficacement que possible.

C'est que les Allemands ont fait, fait-il remarquer. Rejetant les querelles stériles à propos de la science allemande[17] [Nordmann 1916a], il pense qu'on doive les imiter : « le meilleur moyen de nous défendre dans la paix comme dans la guerre contre un ennemi tenace est d'apprendre d'abord à le bien connaître et d'imiter ce qu'il fait de bon » [Nordmann 1917d, 702]. « Allons-nous renoncer à radiographier nos blessés sous prétexte que les rayons X ont été inventés par Rœntgen ? » [Nordmann 1916e, 226]. « Cambrioler » les idées françaises, c'est d'ailleurs ce à quoi les Allemands eux-mêmes sont passés maîtres : « n'eût-il pas mieux valu que ces inventions françaises n'eussent jamais été faites ? Non, si nous savons les généraliser intelligemment jusqu'à dominer l'ennemi non seulement par le cœur, mais aussi par la technique et l'acier[18] » [Nordmann 1916k, 228]. Si les Allemands pervertissent les inventions françaises ou non, c'est parce que leur système moral est inadéquat. Les Français doivent tirer partie de leur supériorité morale et prendre exemple sur les Allemands dans le domaine de la technique. Il faut « louer comme il convient ce qui, dans l'ordre de la technique, est louable chez un ennemi dont on sait assez que, dans l'ordre moral, il a jeté le défi aux lois les plus élémentaires de l'humanité » [Nordmann 1917b, 696].

---

[17] « En dénigrant systématiquement un adversaire qui est toujours debout, on se rabaisse soi-même » [Nordmann 1917b, 695].

L'Allemagne avant la guerre avait déjà su « industrialiser sa science » [Nordmann 1917d, 697]. A Iéna, l'industrie d'optique Zeiss, où « 1500 ouvriers dont le cerveau collectif était constitué par une vingtaine de mathématiciens, physiciens et chimistes » [Nordmann 1917d, 702] doit servir de principe d'organisation industrielle de la science française : « Que voyons-nous à la base de l'œuvre industrielle accomplie à Iéna comme aussi de la plupart des grandes industries allemandes [...] ? *Des recherches de science pure suivies et accompagnées de recherches pratiques de laboratoire* » [Nordmann 1917d, 703, souligné par l'auteur]. A l'instar d'un Henri Le Châtelier, il affirme que les scientifiques doivent être les capitaines de l'industrie. En France, par contraste, la situation dépeinte par Nordmann se rapproche de celle décriée par Claude :

> la science et les savans ne sont pas appréciés comme ils le sont à l'étranger [...]. Si le grand public croit un peu à la science, il n'en est généralement pas de même des pouvoirs publics, ni des chefs d'industrie. Jamais les hommes de sciences ne sont consultés sur les mesures d'intérêt public, même les plus étroitement liées à la science, comme les questions d'organisation de l'enseignement. Dans l'industrie, il en est généralement de même. Dans l'armée même, — j'en sais quelque chose, — beaucoup ont considéré comme anormal et presque scandaleux que certains problèmes militaires d'ordre exclusivement scientifique fussent abordés et, — scandale encore plus grand ! — résolus par des hommes de science [Nordmann 1917d, 704].

Comme on l'a dit, ce n'est pas tant la science dans la guerre qui importe à Nordmann que la guerre scientifique. Puisque la guerre est scientifique, il lui semble normal que les scientifiques soient amenés à intervenir sur des questions de tous ordres, y compris la tactique et la stratégie militaire. Il montre par exemple le peu d'effet des bombardements allemands aériens ou balistiques. Il critique ceux qui croient que le 75 peut remplir tous les besoins de l'artillerie comme ceux qui pensent qu'un seul type d'avion peut servir au bombardement, à la chasse et à la reconnaissance [Nordmann 1918h, 938 ; voir aussi Nordmann 1916h]. « Personnellement convaincu qu'une décision nette ne pourra être obtenue que militairement, et sur notre front » [Nordmann 1916m, 704], il pense parfois que celle-ci ne peut venir que des savants, capable d'inventer une « machine à finir la guerre » [Nordmann 1918h, 935]. Nordmann ne cache rien pourtant les difficultés d'une telle tâche :

---

[18] Le thème des voies françaises non suivies par les Français revient souvent chez Nordmann : « nul n'est prophète en son pays est un proverbe français, trop français même » [Nordmann 1918g, 466]. Pour augmenter la puissance de notre artillerie, écrit-il aussi, « il nous suffira d'entrer dans cette voie si brillamment parcourue par nos Alliés britannique et qui est une voie française, comme Suez, comme Panama… comme tant d'autres chemins tracés dans la matière ou dans l'esprit » [Nordmann 1918g, 469].

> Il est bien difficile de sortir de toutes les contradictions que la guerre traîne dans les plis de sa robe sanglante. Il est bien difficile de savoir si le mieux est ou n'est pas l'ennemi du bien, de savoir si aujourd'hui doit ou non céder le pas à demain.
>
> Et pourtant ces problèmes sont importants : une décimale de plus ou de moins dans le chiffre qui exprime l'efficacité d'un engin peut décider du dort des batailles [...]. C'est ainsi que l'incertitude, l'équilibre si difficile à garder entre le possible et le réel, font que la tactique, la technique de guerre qu'on pourrait croire la plus positive des sciences, finissent par devenir de petites choses presque aussi conjecturales [...] que la métaphysique [Nordmann 1918h, 945].

Mais si la complexité de la guerre scientifique peut faire douter de la fécondité de la méthode scientifique, il est clair que c'est là que réside le seul espoir. Dès 1916, Nordmann insiste sur l'importance d'une vision de scientifique afin de penser la France d'après-guerre. Pour développer sa richesse future, il faut augmenter les rendements agricoles grâce à la mécanisation [Nordmann 1916c ; 1916d]. Il faut prêter attention au développement économique de la nation, condition essentielle de l'épanouissement des arts, des lettres et des sciences [Nordmann 1917e]. Parce qu'il veut croire que la mobilisation scientifique aura *in fine* des conséquences bénéfiques pour l'humanité, Nordmann suivra avec grand intérêt, après guerre, les tentatives de reconstruction des associations internationales, consacrant plusieurs articles à la Croix Rouge[19], à la réforme du calendrier[20] et au Conseil international de recherche (CIR) [Nordmann 1919k].

**Conclusion**

Les historiens de la Grande Guerre nous avaient montré que si l'idée de progrès est d'abord remise en cause en 1914, si on déplore l'instrumentalisation de la science, si on fustige les barbaries de la science allemande, la réaction fait long feu, même au cœur de la RDM qui avait si longtemps résisté à la prétention scientifique de vouloir fonder une nouvelle morale. De fait, la rénovation des structures scientifiques (l'Académie, les commissions et le ministère des inventions, les nouveaux instituts et laboratoires) semble faire l'unanimité. Pourtant, ni la religion scientiste ni l'utilitarisme pragmatique à outrance ne s'impose [Prochasson & Rasmussen 1996].

La trajectoire de Nordmann permet d'éclairer ces processus et montre à quel point les valeurs conservatrices et idéalistes de l'avant-guerre ont pu servir de moule à une nouvelle conception de la science. En étudiant la manière dont il acclimate un discours relativement conservateur aux conditions nouvelles suscitées par la guerre, on comprend mieux que l'opinion puisse être relativement unanime. Quelque chose de neuf est né : une opinion qui voit dans la science une entité amorale, parfois dangereuse, mais qu'on ne peut guère négliger qu'au péril de la nation. Les valeurs positives de la science sont réinvesties dans la nation qui devient garante des dérives toujours possibles. Plus tard, les organisations internationales seront réinvesties de ces valeurs qui encadrent une science qui, bien que dangereuse, n'a rien perdu de son attrait.

---

[19] La maladie : « Voilà un nouveau champ de bataille où le devoir est de combattre maintenant » [Nordmann 1919e, 458]. Voir ses articles sur la grippe et la tuberculose [1918i ; 1918j ; 1918k].
[20] Il présente d'ailleurs l'Association géodésique internationale fondée en 1866 comme « une cellule embryonnaire de la Société des Nations, une cellule née de la science » [Nordmann 1919a, 461]. Voir aussi [Nordmann 1919b]

C'est en novembre 1919 que la guerre scientifique prend fin pour Nordmann. Envoyé à Londres pour suivre l'assemblée du CIR, il découvre avec stupéfaction les nouvelles recherches atomiques. Il tourne alors la page et retourne à ses anciennes amours :

> Depuis cinq ans, la plupart des hommes de sciences avaient quitté les hautes régions de la recherche spéculative pour se consacrer, — chacun suivant ses forces, — à la défense de la patrie. Ici même depuis lors, et parce que nous tâchons de suivre fidèlement le mouvement des idées, ce sont des questions concrètes, des problèmes de technique pratique et de science appliquée qui ont fourni exclusivement les sujets de nos chroniques. […] Et il fallait quelqu'effort d'imagination pour se souvenir qu'elles n'étaient après tout que des reflets, — transformés sur le dur miroir des nécessités vitales, — de ces hautes études qui s'appellent la mécanique rationnelle, l'acoustique et l'optique, la chimie pure, la physiologie, la biologie.

> Maintenant les hommes de sciences regagnent leur laboratoire familier. Ou du moins, car beaucoup d'entre eux ne l'avaient jamais quitté, ils s'apprêtent à refermer doucement la fenêtre qu'ils y avaient ouverte un moment et qui leur donnait vue sur le champ de bataille. Les tours d'ivoire de la science pure naguère dépeuplées, et abandonnées à leur solitude au milieu des nuages, voient revenir ceux qui les avaient désertées, un peu plus courbés peut-être sous le fardeau de cinq années d'angoisse et avec bien des vides dans leurs rangs ; mais avec toujours au cœur cette flamme ardente de la curiosité, avec toujours cette passion de savoir pour savoir, et non point seulement pour pouvoir [Nordmann 1919j, 214].

La parenthèse semble s'être refermée. Pourtant une analyse approfondie des écrits de Nordmann dans les années 1920 montrerait qu'il ne désarme pas. Pendant la visité controversée d'Einstein à Paris en 1922, Nordmann sera encore une fois sur la ligne de front, non plus pour combattre les Allemands, mais les préjugés et l'ignorance[21]. Sans renier Poincaré, Nordmann n'aura cependant plus la prétention de vouloir fonder un quelconque sentiment religieux sur les sciences :

> On ne trouvera pas, écrit-il dans un livre publié en 1923, les habituelles niaiseries et fadaises sur les beaux sentiments que doit conférer l'étude du ciel. Je n'y crois pas ne l'ayant guère observé autour de moi ; et il y a longtemps qu'Henri Poincaré a justement proclamé la séparation de la Science et de la Morale. […] Nous examinerons donc, avec sincérité, si cette harmonie immense du firmament stellaire doit seulement éblouir nos regards et si elle ne peut pas aussi apaiser nos cœurs. Sans inutile désespoir, sans illusions qui enivrent, nous chercherons le divin sous les formes sensibles [Nordmann 1923a, 6 & 11].

## Bibliographie

---

[21] [Biezunski 1991]. Voir aussi les articles de Nordmann dans l'*Illustration* (28 mai 1921 et 15 avril 1922) et dans le *Matin*, et les deux ouvrages très populaires qu'il publie à ce sujet : [Nordmann 1921a & 1924a].

NORDMANN, Charles (1918j) : La lutte sociale contre la tuberculose, RDM, 47, **xxx**, 1918.

NORDMANN, Charles (1918k) : La grippe, RDM, 48, **xxx**, 1918.

NORDMANN, Charles (1918l) : Une révolution en chirurgie, RDM, 48, **xxx**, 1918.

NORDMANN, Charles (1919a) : Le congrès de la paix et le calendrier, RDM, 49, 459-468, 1919.

NORDMANN, Charles (1919b) : Pour une réforme chronologique, RDM, 49, **xxx**, 1919.

NORDMANN, Charles (1919c) : Le relevage des navires torpillés, RDM, 50, 459-470, 1919.

NORDMANN, Charles (1919d) : Un modèle d'organisation : le secours de guerre, RDM, 50, 935-, 1919.

NORDMANN, Charles (1919e) : Une croisade contre la maladie : le congrès de Cannes et le bureau d'hygiène mondial, RDM, 51, 454-465, 1919.

NORDMANN, Charles (1919f) : Les progrès de la T.S.F. et la guerre, RDM, 51, 927-938, 1919.

NORDMANN, Charles (1919g) : Quelques progrès guerriers de l'aéronautique et de l'aérologie, RDM, 52, 693-702, 1919.

NORDMANN, Charles (1919h) : Ecoute sous-marine ; écoute souterraine, RDM, 52, 935-945, 1919.

NORDMANN, Charles (1919i) : Les vitamines, RDM, 53, 686-697, 1919.

NORDMANN, Charles (1919j) : Au royaume de l'infiniment petit, RDM, 54, 214-225, 1919.

NORDMANN, Charles (1919k) : Impressions anglaises, RDM, 54, 458-468, 1919.

NORDMANN, Charles (1919l) : L'électron, RDM, 54, 935-946, 1919.

NORDMANN, Charles (1920a) : Les origines du repérage par le son, Revue scientifique 58, 737-740, 1920.

NORDMANN, Charles (1920b) : Les distances des étoiles, RDM, 55, 457-468, 1920.

NORDMANN, Charles (1920c) : Amas stellaires et nébuleuses, RDM, 55, 929-940, 1920.

NORDMANN, Charles (1920d) : La physiologie des étoiles, RDM, 56, 457-468, 1920.

NORDMANN, Charles (1920e) : Le rôle de la sécrétion interne, RDM, 56, 931-942, 1920.

NORDMANN, Charles (1920f) : L'énigme martienne, RDM, 57, 453-464, 1920.

NORDMANN, Charles (1920g) : Les messages de l'infini, RDM, 57, 914-925, 1920.

NORDMANN, Charles (1920h) : Le Soleil et l'aimant terrestre, RDM, 58, **xxx**, 1920.

NORDMANN, Charles (1920i) : L'action électrique du soleil, RDM, 58, **xxx**, 1920.

NORDMANN, Charles (1920j) : La nouvelle étoile du Cygne, RDM, 59, **xxx**, 1920.

NORDMANN, Charles (1920k) : Fléaux exotiques, RDM, 59, **xxx**, 1920.

NORDMANN, Charles (1920l) : Transports aériens, RDM, 60, **xxx**, 1920.

NORDMANN, Charles (1921?) : dans l'*Illustration* (28 mai 1921 et 15 avril 1922) et dans le *Matin*.

NORDMANN, Charles (1921a) : *Einstein et l'univers : une lueur dans le mystère des choses*, Hachette, 1921.

NORDMANN, Charles (1921b) : L'étude des fluctuations lumineuses de étoiles, RDM, 61, **xxx**, 1921.

NORDMANN, Charles (1921c) : Impressions de Roumanie, RDM, 61, **xxx**, 1921.

NORDMANN, Charles (1921d) : Le phare à éclipse de l'infini, RDM, 62, **xxx**, 1921.

NORDMANN, Charles (1921e) : Qu'est-ce que les rayons X ?, RDM, 62, **xxx**, 1921.

NORDMANN, Charles (1921f) : Le compas lumineux, RDM, 63, **xxx**, 1921.

NORDMANN, Charles (1921g) : Nouveaux horizons en médecine, RDM, 63, **xxx**, 1921.

NORDMANN, Charles (1921h) : A propos de boxe, RDM, 64, **xxx**, 1921.

NORDMANN, Charles (1921i) : Pour préluder à l'étude d'Einstein, RDM, 64, **xxx**, 1921.

NORDMANN, Charles (1921j) : L'espace et le temps selon Einstein, RDM, 65, 313-343, 1921.

NORDMANN, Charles (1921k) : La mécanique d'Einstein, RDM, 65, 925-946, 1921.